\documentclass[12pt,a4paper]{article}

\usepackage{indentfirst}

\usepackage{amsmath,amssymb,amsfonts,amsthm,graphics}
\usepackage{makeidx}
\begin{document}

\title{\bf \Large Two rainbow connection numbers\\ and the parameter
$\sigma_k(G)$\footnote{Supported by NSFC.} }
\author{\small Jiuying Dong, Xueliang Li\\
\small Center for Combinatorics and LPMC-TJKLC \\
\small Nankai University, Tianjin 300071, China\\
\small Email: jiuyingdong@126.com; lxl@nankai.edu.cn}

\date{}

\maketitle{}

\begin{abstract}

The rainbow connection number $rc(G)$ and the rainbow
vertex-connection number $rvc(G)$ of a graph $G$ were introduced by
Chartrand et al. and Krivelevich and Yuster, respectively. Good
upper bounds in terms of minimum degree $\delta$ were reported by
Chandran et al., Krivelevich and Yuster, and Li and Shi. However, if
a graph has a small minimum degree $\delta$ and a large number of
vertices $n$, these upper bounds are very large, linear in $n$.
Hence, one may think to look for a good parameter to replace
$\delta$ and decrease the upper bounds significantly. Such a natural
parameter is $\sigma_k$. In this paper, for the rainbow connection
number we prove that if $G$ is a connected graph of order $n$ with
$k$ independent vertices, then $rc(G)\leq
3k\frac{n-2}{\sigma_k+k}+6k-4$. For the rainbow vertex-connection
number, we prove that $rvc(G)\leq \frac{(4k+2k^{2})n}{\sigma_k+k}+5k
$ if $\sigma_k\leq 7k$ and $\sigma_k\geq 8k$, and $rvc(G)\leq
\frac{(\frac{38k}{9}+2k^{2})n}{\sigma_k+k}+5k$ if $7k<\sigma_k< 8k$.
Examples are given showing that our bounds are much better than the
existing ones, i.e., for the examples $\delta$ is very small but
$\sigma_k$ is very large, and the bounds are $rc(G)< 9k-3$ and
$rvc(G)\leq 9k+2k^{2}$ or $rvc(G)\leq\frac{83k}{9}+2k^{2}$, which
imply that both $rc(G)$ and $rvc(G)$ can be upper bounded by
constants from our upper bounds, but linear in $n$ from the
existing ones.\\[3mm]
\end{abstract}

\section{ Introduction}

All graphs under our consideration are finite, undirected and
simple. For notations and terminology not defined here, we refer to
[2]. Let $G$ be a graph. The length of a path $P$ in $G$ is the
number of edges of $P$. The distance between two vertices $u$ and
$v$ in $G$, denoted by $d(u, v)$, is the length of a shortest path
connecting them in $G$. If there is no path connecting $u$ and $v$,
we set $d(x,y):=\infty$. For two subsets $X$ and $Y$ of $V$, an
$(X,Y)$-path is a path which connects a vertex of $X$ and a vertex
of $Y$, and whose internal vertices belong to neither $X$ nor $Y$.
We use $E[X,Y]$ to denote the set of edges of $G$ with one end in
$X$ and the other end in $Y$, and $e(X,Y)=|E[X,Y]|$.

Let $c: E(G)\rightarrow \{1,2,\cdots,k\}, k\in N$ be an
edge-coloring, where adjacent edges may be colored the same. A graph
$G$ is rainbow edge-connected if for every pair of distinct vertices
$u$ and $v$ of $G$, $G$ has a $u-v$ path $P$ whose edges are colored
with distinct colors. This concept was introduced by Chartrand et
al. [5]. The minimum number of colors required to rainbow color a
connected graph is called its rainbow connection number, denoted by
$rc(G)$. Observe that if $G$ has $n$ vertices then $rc(G)\leq n-1$.
Clearly, $rc(G) \geq diam(G)$, the diameter of $G$. In [5],
Chartrand et al. determined the rainbow connection numbers of
wheels, complete graphs and all complete multipartite graphs. In
[3], Caro et al. got the following theorems and made the following
conjectures.\\

\noindent{\bf Theorem 1 [3].} If $G$ is a connected graph with $n$
vertices and minimum degree $\delta\geq 3$, then $rc(G)<
\frac{5}{6}n$.\\

\noindent{\bf Theorem 2 [3].} If $G$ is a connected graph with $n$
vertices and minimum degree $\delta$, then $rc(G)\leq
min\{\frac{ln\delta}{\delta}n(1+o_\delta(1)),
n\frac{4ln\delta+3}{\delta}\}$.\\

\noindent{\bf Conjecture 1 [3].} If $G$ is a connected graph with
$n$ vertices and minimum degree $\delta\geq 3$, then $rc(G)<
\frac{3}{4}n$.\\

In [6], Krivelevich and Yuster got the following upper bound, which
looks much simpler than Theorem 2.\\

\noindent{\bf Theorem 3 [6].}  A connected graph $G$ with $n$
vertices has $rc(G)\leq \frac{20n}{\delta}$.\\

In [8], Schiermeyer proved Conjecture 1 and posed the following
challenging problem.\\

\noindent{\bf Problem 1 [8].}  For every $k\geq 2$ find a minimal
constant $c_k$ with $ 0< c_k\leq 1$ such that $rc(G)\leq nc_k$ for
all graphs $G$ with minimum degree $\delta \geq k$. Is it true that
$c_k=\frac{3}{k+1}$ for all $k\geq 27$ ?\\

In [4], Chandran et al. settled this problem, which improves the
result of Theorem 3. They obtained the following result.\\

\noindent{\bf Theorem 4 [4].} For every connected graph $G$ of order
$n$ and minimum degree $\delta$, we have $rc(G)\leq 3n/(\delta + 1)
+ 3$. Moreover, the bound is seen to be tight up to additive factors
by examples given in [3].\\

Before the proof of Theorem 4, they first proved the following
result.\\

\noindent{\bf Theorem 5 [4].} If $D$ is a connected two-way two-step
dominating set in a graph $G$, then $rc(G) \leq rc(G[D]) +
6$.\\

However, if a graph has a small minimum degree $\delta$ but a large
order $n$, then the upper bound $3n/(\delta + 1) + 3$ will be very
large, linear in $n$. But, actually $rc(G)$ could be much less than
the bound, a constant, for which we will give examples later. Hence,
one may think to look for a better parameter to replace $\delta$.
Such a natural parameter is $\sigma_k$, which is defined by
$\sigma_k(G)=min \{d(u_1)+d(u_2)+\cdots+d(u_k) | \
u_1,u_2,\ldots,u_k\in V(G), u_iu_j\not\in E(G),i\neq j,
i,j\in\{1,\cdots,k\}\}$, or simply denoted by $\sigma_k$. Observe
that $\sigma_k$ is monotonically increasing in $k$. So $\sigma_k$
could be very large, which would decrease the upper bound
dramatically. In this paper, we will employ the parameter $\sigma_k$
to get the
following result.\\

\noindent{\bf Theorem 6.} If $G$ is a connected graph of order $n$
with $k$ independent vertices, then $rc(G)\leq
3k\frac{n-2}{\sigma_k+k}+6k-4$.\\

From the following examples, one can see that $\sigma_k$ really
works very well in decreasing the upper bound of $rc(G)$. First of
all, we denote by $K^{*}_{a,b}$ the graph obtained from the complete
bipartite graph $K_{a,b}$ by joining every pair of vertices
in the $b$-part by a new edge.\\

\noindent{\bf Example 1:} Let $\frac{n-2}{k-1}$ be an integer and
let $H=K^{*}_{2, \frac{n-2}{k-1}-2}$,
$H_1=K^{*}_{2,\frac{n-2}{k-1}-1}$, and $H_k=K_1$ with
$V(K_1)=\{v\}$. Take $k-2$ copies of $H$, denoted by $H_2,\cdots,
H_{k-1}$. Label the two non-adjacent vertices of $H_i$ by $x_{i,1},
x_{i,2}$, for $i\in \{1,\cdots,k-1\}$. Now, connect $x_{i,2}$ and
$x_{i+1,1}$ with an edge for every $i\in \{1,\cdots,k-1\}$, and
connect $v$ and $x_{k-1,2}$ with an edge. The resulting graph is
denoted by $G$. From the construction, it is not difficult to check
that for every $v\in V(H_i),i\in\{2,\cdots,k-1\}$, we have
$d(v)=\frac{n-2}{k-1}-1$. In addition,
$d(x_{1,1})=\frac{n-2}{k-1}-1$,
$\sigma_k=(\frac{n-2}{k-1}-1)(k-1)+1=n-k$, and $\delta(G)=1$. From
these facts, one can see that the upper bound of Theorem 4 is
$rc(G)\leq 3n/2 + 3$, which is linear in $n$, nevertheless, the
upper bound in our Theorem 6 is $rc(G) < 9k-4$, which is a constant
when $k$ is small, say 2, 3, etc. Notice that here we can make
$\delta$ be 2,3, etc, simply by adding a few edges properly.\\

\noindent{\bf Example 2:} Let $\frac{\sigma_k}{k}$ be an integer and
let $H=K^{*}_{2, \frac{\sigma_k}{k}-1}$,
$H'=K^{*}_{2,\frac{\sigma_k}{k}}$. Take $t$ copies of $H$, denoted
by $H_1,\cdots, H_t$, and take two copies of $H'$, denoted by $H_0,
H_{t+1}$. Label the two non-adjacent vertices of $H_i$ by $x_{i,1},
x_{i,2}$, for $i\in\{0,1,\cdots,t+1\}$. Now, connect $x_{i,2}$ and
$x_{i+1,1}$ for $i\in \{0,\cdots,t+1 \}$ with an edge. The resulting
graph $G$ has $n=(t+2)(\frac{\sigma_k}{k} +1)+2$ vertices. It is
straightforward to verify that for $i\in\{1,2,\cdots,t\}$ and any
$v\in V(H_i)$, we have $d(v)=\frac{\sigma_k}{k}$. In addition,
$d(x_{0,1})=d(x_{t+1,2})=\frac{\sigma_k}{k}$, and
$diam(G)=d(x_{0,1},x_{t+1,2})=3t+5$. From
$3k\frac{n-2}{\sigma_k+k}-1=3t+5$, and $rc(G)\geq diam(G)$, one can
see that the bound $rc(G)\leq 3k\frac{n-2}{\sigma_k+k}+6k-4$ of
Theorem 6 could be seen to be tight up to additive factors $6k-3$
when $k$ is small.\\

Let $c: V(G)\rightarrow \{1,2,\cdots,k\}, k\in N$ be a
vertex-coloring, where adjacent vertices may be colored the same. A
graph $G$ is rainbow vertex-connected if for every pair of distinct
vertices $u$ and $v$ of $G$, $G$ has a $u-v$ path $P$ whose internal
vertices are colored with distinct colors. The minimum number of
colors required to rainbow color a connected graph is called the
rainbow vertex-connection number of $G$, denoted by $rvc(G)$. The
concept of rainbow vertex-connection number was introduced by
Krivelevich and Yuster [6]. It is obvious that $rvc(G)\leq n-2$ and
$rvc(G) \geq diam(G)-1$. In [6], Krivelevich
and Yuster obtained the following result:\\

\noindent{\bf Theorem 7 [6].} A connected graph $G$ of order $n$
with minimum degree $\delta$ has $rvc(G) <\frac{11n}
{\delta}$.\\

In [7], Li and Shi improved the above bound and got the following
result:\\

\noindent{\bf Theorem 8 [7].} A connected graph $G$ of order $n$
with minimum degree $\delta$ has $rvc(G) \leq
\frac{4n}{\delta+1}+C(\delta)$ for $\delta\geq 6$, where $C(\delta)=
e^{\frac{3log(\delta^{3}+2\delta^{2}+3)-3(log3-1)}{\delta-3}}-2$.
And $rvc(G) \leq \frac{3n}{4}-2$ for $\delta =3 $, $rvc(G) \leq
\frac{3n}{5}-\frac{8}{5}$ for $\delta =4 $, $rvc(G) \leq
\frac{n}{2}-2$ for $\delta =5 $.\\

Similar to the edge-coloring case, if we use the parameter
$\sigma_k(G)$ to replace $\delta $, the upper bound of $rvc(G)$ can
also be dramatically improved, see the following result:\\

\noindent{\bf Theorem 9.} Let $G$ is a connected graph of order $n$
with $k$ independent vertices. Then $rvc(G)\leq
\frac{(4k+2k^{2})n}{\sigma_k+k}+5k $ if $\sigma_k\leq 7k$ and
$\sigma_k\geq 8k$; whereas $rvc(G)\leq
\frac{(\frac{38k}{9}+2k^{2})n}{\sigma_k+k}+5k$ if $7k<\sigma_k<
8k$.\\

From Example 1 one can see that there are infinitely many graphs $G$
satisfying $\sigma_k\geq n-k$ and $\delta$ is small, which means
that $rvc(G)\leq 9k+2k^{2}$ or $rvc(G)\leq \frac{83}{9}k+2k^{2}$,
which are constants, however the bounds in Theorems 7 and 8 give us
$rvc(G) \leq 11n/\delta$ and $rvc(G) \leq
\frac{4n}{\delta+1}+C(\delta)$, which are linear in $n$.\\

The rainbow connection numbers have applicable background. They can
be used in secure transfer of classified information between
agencies. Suppose we have a communication network $G$, and we want
to transfer information between any two agencies along a route in
the network in such a way that each link on the route is assigned a
distinct channel. The aim is to use as few distinct channels as
possible in our network. The question is what is the minimal number
of channels we have to use. $rc(G)$ is just the number we want.\\

The following notions are needed in what follows, which could be
found in [4, 6]. Given a graph $G$, a set $D\subseteq V(G)$ is
called a $k$-step dominating set of $G$, if every vertex in $G$ is
at a distance at most $k$ from $D$. Further, if $D$ induces a
connected subgraph of $G$, it is called a connected $k$-step
dominating set of $G$. The $k$-step open neighborhood of a set
$D\subseteq V(G)$ is $N^{k}(D):=\{x\in V(G)| d(x,D)=k \}$,
$k=\{0,1,2,\cdots\}$. A connected two-step dominating set $D$ in a
graph $G$ is called a connected two-way two-step dominating set if
every pendant vertex of $G$ is included in $D$ and every vertex in
$N^{2}(D)$ has at least two neighbors in $N^{1}(D)$. We call a
two-step dominating set $k$-strong if every vertex in $N^{2}(D)$ has
at least $k$ neighbors in $N^{1}(D)$.

\section{ Proof of Theorem 6 }

\noindent{\bf Theorem 6.} If $G$ is a connected graph of order $n$
with $k$ independent vertices,  then $rc(G)<
3k\frac{n-1}{\sigma_k+k}+6k-3$.\\

\noindent{\bf Proof.} As $rc(G)\leq n-1$, if $\sigma_k\leq 2k$, then
$3k\frac{n-1}{\sigma_k+k}+6k-3\geq n+6k-5\geq n+7 (k\geq 2)$, the
theorem is true. So we may assume that $\sigma_k\geq 2k+1$. First we
see the following two claims: \\

\noindent{\bf Claim 1.} $G$ has a connected two-step dominating set
$D$ such that $|D|\leq
3k\frac{n-|N^{2}(D)|-1}{\sigma_k+k}+3k-5$.\\

\noindent{\bf Proof.} Let $u_1,u_2,\cdots,u_k $ be independent
vertices of $G$ and $d(u_1)\geq d(u_2)\geq \cdots\geq d(u_k)$. Say
$D=\{u_1\}$, we have $|D\cup N^{1}(D)| \geq
\lceil\frac{\sigma_k}{k}\rceil+1$.\\

\noindent{\bf Case 1.1.} $\alpha (G[N^{3}(D)])\geq k.$

Let $v_1,v_2,\cdots,v_k$ be independent vertices of $G[N^{3}(D)]$
and $d(v_1)\geq d(v_2)\geq\cdots\geq d(v_k)$. Since $d(v_1,D)=3$,
let $P=v_1v^{2}_1v^{1}_1v^{0}_1$ be a shortest $v_1-D$ path where
$v^{2}_1\in N^{2}(D), v^{1}_1\in N^{1}(D), v^{0}_1\in D$. (Latter,
we will omit this note). Say $D=\{u_1,v_1,v^{2}_1,v^{1}_1\}$. When
the vertex $v_1$ was put to $D$, $|D\cup N^{1}(D)|$ increases by at
least $\lceil\frac{\sigma_k}{k}\rceil+1$. If $\alpha
(G[N^{3}(D)])\geq k $, we continue the above manipulation. When
$\alpha (G[N^{3}(D)])\leq k-1 $, we may get $|D| \leq 3(\frac{|D\cup
N^{1}(D)|}{\lceil\frac{\sigma_k}{k}\rceil+1}-1)+1\leq
3k\frac{n-|N^{2}(D)|-|N^{3}(D)|}{{\sigma_k}+k}-2$. Let
$\{x_1,x_2,\cdots,x_t\}$ be a maximum independent set of
$G[N^{3}(D)]$. $P_1=x_1x^{2}_1x^{1}_1x^{0}_1$ be a shortest $x_1-D$
path, say $D=D\cup\{x_1,x^{2}_1,x^{1}_1\}$. Similarly, let
$P_2=x_2x^{2}_2x^{1}_2x^{0}_2$ be a shortest $x_2-D$ path, say
$D=D\cup\{x_2,x^{2}_2,x^{1}_2\},\cdots.$ Finally, let
$P_t=x_tx^{2}_tx^{1}_tx^{0}_t$ be a shortest $x_t-D$ path, say
$D=D\cup\{x_t,x^{2}_t,x^{1}_t\} $. We may see that $D$ is a
connected two-step dominating set and $|D| \leq
3k\frac{n-|N^{2}(D)|-|N^{3}(D)|}{{\sigma_k}+k}-2 +3(k-1)\leq
3k\frac{n-|N^{2}(D)|-1}{\sigma_k+k}+3k-5$. So the claim is
true.\\

\noindent{\bf Case 1.2.} $\alpha (G[N^{3}(D)])\leq k-1.$

Let $\{x_1,x_2,\cdots,x_t\}$ be a maximum independent set of
$G[N^{3}(D)]$. Similar to the proof of the latter part of Case 1.1,
we get $|D| \leq 1+3t\leq1+3(k-1)$. It is obvious that the claim is
true again. \qed \\

\noindent{\bf Claim 2.} If $\sigma_k\geq 2k+1$, then $G$ has a
connected two-way two-step dominating set $D$ such that $|D|\leq
3k\frac{n-2}{\sigma_k+k}+6k-9$.\\

We look at the connected two-step dominating set $D$ of Claim 1. As
$\sigma_k\geq 2k+1$, $N^{1}(D)$ has at most $k-1$ pendant vertices.
We put the $k-1$ pendant vertices to $D$. So $|D|\leq
3k\frac{n-|N^{2}(D)|-1}{\sigma_k+k}+3k-5+k-1=
3k\frac{n-|N^{2}(D)|-1}{\sigma_k+k}+4k-6$. Note that $N^{1}(D)$ has
no pendant vertices. If for each vertex $v\in N^{2}(D)$, $e(v,
N^{1}(D))\geq 2$, then $D$ is exactly the required connected two-way
two-step dominating set $D$, and the claim is true. Therefore, we
may assume that $\exists v\in N^{2}(D), e(v, N^{1}(D))=1$.\\

\noindent{\bf Case 2.1.} There exists an independent set
$\{v_1,v_2,\cdots,v_k\}$ in $G[N^{2}(D)]$ such that $ e(v_1,
N^{1}(D))= e(v_2, N^{1}(D))=\cdots = e(v_k, N^{1}(D))=1$.

Suppose that $d(v_1)\geq d(v_2)\cdots\geq d(v_k)$. As $d(v_1,D)=2$,
let $P=v_1v^{1}_1v^{0}_1$ be a shortest $v_1-D$ path, where $v_1$
has at least $\lceil\frac{\sigma_k}{k}\rceil-1$ neighbors in
$N^{2}(D)$. When we put vertex $v_1$ to $D$, $| N^{2}(D)|$ reduces
by at least $\lceil\frac{\sigma_k}{k}\rceil$. If $G[N^{2}(D)]$ still
has an independent set $\{v_1,v_2,\cdots,v_k\}$ such that $ e(v_1,
N^{1}(D))= e(v_2, N^{1}(D))=\cdots = e(v_k, N^{1}(D))=1$, we
continue the above manipulation, until $G[N^{2}(D)]$ has no
independent set $\{v_1,v_2,\cdots,v_k\}$ such that $ e(v_1,
N^{1}(D))= e(v_2, N^{1}(D))=\cdots = e(v_k, N^{1}(D))=1$. Thus $D$
increases by at most $2|N^{2}(D)|/\lceil\frac{\sigma_k}{k}\rceil
\leq\frac{2k|N^{2}(D)|}{\sigma_k}$. Hence $|D|\leq
3k\frac{n-|N^{2}(D)|-1}{\sigma_k+k}+4k-6+\frac{2k|N^{2}(D)|}{\sigma_k}<
3k\frac{n-1}{\sigma_k+k}+4k-6$. Here, $ N^{2}(D)$ can be partitioned
into two parts $N_1^{2}(D)$ and $N_2^{2}(D)$, for  $\forall v\in
N_1^{2}(D), e(v, N^{1}(D))\geq 2$, and $\forall v\in N_2^{2}(D),
e(v, N^{1}(D))=1$ and $\alpha (G[N_2^{2}(D)])\leq k-1$, where
$|N_1^{2}(D)|\geq 0, | N_2^{2}(D)|\geq 0$. In the same way as
before, we can arrive at that $D$ is a connected two-way two-step
dominating set such that $|D|
<3k\frac{n-1}{\sigma_k+k}+4k-6+2(k-1)=3k\frac{n-1}{\sigma_k+k}+6k-8$.
So the claim is true.\\

\noindent{\bf Case 2.2.} There does not exist any independent set
$\{v_1,v_2,\cdots,v_k\}$ in $G[N^{2}(D)]$ such that $ e(v_1,
N^{1}(D))= e(v_2, N^{1}(D))=\cdots = e(v_k, N^{1}(D))=1$.

We partition $ N^{2}(D)$ into two parts $N_1^{2}(D)$ and
$N_2^{2}(D)$, for $\forall v\in N_1^{2}(D), e(v, N^{1}(D))\geq 2$,
and $\forall v\in N_2^{2}(D), e(v, N^{1}(D))=1$ and $\alpha
(G[N_2^{2}(D)])\leq k-1$, where $| N_1^{2}(D)|\geq 0, |
N_2^{2}(D)|\geq 0$. Similarly, we can get a connected two-way
two-step dominating set such that $|D|< 3k\frac{n-
|N^{1}(D)|-1}{\sigma_k+k}
+4k-6+2(k-1)<3k\frac{n-2}{\sigma_k+k}+6k-8$. The claim is again
true. \qed \\

Observe that the connected two-way two-step dominating set $D$ can
be rainbow colored, using $|D|-1$ colors by ensuring that every edge
of some spanning tree gets distinct colors. According to Claim 2 and
Theorem 5, the upper bound follows immediately. \qed

\section{ Proof of Theorem 9}

We first recall the following Lemma 1 and prove Lemma 2, as we need
them in the proof of our theorem.\\

\noindent{\bf Lemma 1 (The Lov\'{a}sz Local Lemma [1]).} Let
$A_1,A_2, \cdots, A_n$ be the events in an arbitrary probability
space. Suppose that each event $A_i$ is mutually independent of a
set of all the other events $A_j$ but at most $d$, and that $P[A_i]
\leq p$ for all $1\leq i\leq n$. If $ep(d + 1) <1$, then
$Pr[\bigwedge^{n}_{i=1}\overline{A_i}]> 0$.\\

\noindent{\bf  Lemma 2.} If $G$ is a connected graph of order $n$
with $k$ independent vertices, then $G$ has a connected spanning
subgraph $G'$ which has the same value of $\sigma_k$ as $G$ and $
e(G')< n(\sigma_k-k+1)+\frac{kn}{\sigma_k+k}$.\\

\noindent{\bf Proof.}  For convenience, we denote by $I_k$ an
independent set $\{v_1,v_2,\cdots,v_k\}$ that satisfies
$d(v_1)+d(v_2)+\cdots d(v_k)=\sigma_k$. We delete the edges of $G$
as soon as possible and get $H$ such that $\sigma_k(H)=\sigma_k$,
but for each edge $e\in E(H), \sigma_k(H-e)<\sigma_k$. Hence each
edge of $H$ is incident to some vertex of some $I_k$. Suppose that
$H$ has mutually disjoint independent sets
$I^{1}_{k},I^{2}_{k},\cdots,I^{a}_{k}$, that is, $I^{i}_{k}\cap
I^{j}_{k}=\phi, i,j \in\{1,2,\cdots,a\}$. Let $H'=
H-\bigcup_{i=1}^{a}I^{i }_{k}$. Then for each edge $e\in E(H')$, at
least one of its ends $w$ is in some $I_k$ and $d_w(H)\leq
\sigma_k-(k-1)$. Suppose that all edges of $E(H')$ are incident to
$b$ vertices each of which is in some $I_k$. We know that for each
vertex $v\in V(H'\cap I_k), d(v)\leq \sigma_k-(k-1)$. Note that
$b\leq n-ka$, so $e(H)<a\sigma_k+b(\sigma_k-k+1)\leq
a\sigma_k+(n-ka)(\sigma_k-k+1)=a(1-k)(\sigma_k-k)+n(\sigma_k-k+1)$.
If $H$ has $t$ connected components $H_1,H_2,\cdots, H_t$, then for
any $k$ vertices $v_{i_1},v_{i_2},\cdots,v_{i_k}$, each $v_{i_j}$
taken from the corresponding component $H_{i_j}$. We have
$d(v_{i_1})+d(v_{i_2})+\cdots +d(v_{i_k})\geq \sigma_k$. So
$|V(H_{i_1})|+|V(H_{i_2})|+\cdots+|V(H_{i_k})|\geq \sigma_k+k$.
Thus, $t\leq \lceil \frac{kn}{\sigma_k+k}\rceil$. That is, $H$ has
at most $\lceil\frac{kn}{\sigma_k+k}\rceil$ connected components.
Therefore, we get $e(G')< a(1-k)(\sigma_k-k)+n(\sigma_k-k+1) +
\lceil\frac{kn}{\sigma_k+k}\rceil-1 <
n(\sigma_k-k+1) +\frac{kn}{\sigma_k+k}$. The claim follows. \qed \\

\noindent{\bf The proof of Theorem 9.}  Since $G$ is a connected
graph of order $n$, we know $rvc(G)\leq n-2$. As
$4k+2k^{2}-8k=2k(k-2)\geq 0$, if $\sigma_k\leq 7k$, the result is
obvious. So we may assume $\sigma_k\geq 7k+1$.\\

\noindent{\bf Claim 3.} $G$ has a connected
$\lceil\frac{\sigma_k}{2k}\rceil$-strong two-step dominating set $D$
such that $|D|< 4k\frac{n-1}{\sigma_k+k}+5k-6$.\\

\noindent{\bf Proof.} We look at the set $D$ in Claim 1 of the proof
for Theorem 6. If for each vertex $v\in N^{2}(D)$, $e(v,
N^{1}(D))\geq \lceil\frac{\sigma_k}{2k}\rceil$, then $D$ is exactly
the required dominating set. So we assume that there exists a $v\in
N^{2}(D)$ such that $e(v, N^{1}(D))\leq
\lceil\frac{\sigma_k}{2k}\rceil-1$.\\

\noindent{\bf Case 3.1.}  There exists an independent set
$\{v_1,v_2,\cdots,v_k\}$ in $G[N^{2}(D)]$ such that $ e(v_i,
N^{1}(D)) \leq \lceil\frac{\sigma_k}{2k}\rceil-1, i\in \{1,2,\cdots,
k\}$.

We assume $d(v_1)\geq d(v_2)\geq\cdots\geq d(v_k)$. Let
$P=v_1v^{1}_1v^{0}_1$ be a shortest $v_1-D$ path, where $v_1$ has at
least
$\lceil\frac{\sigma_k}{k}\rceil-\lceil\frac{\sigma_k}{2k}\rceil+1
\geq\lfloor\frac{\sigma_k}{2k}\rfloor+1$ neighbors in $N^{2}(D)$. When
we put vertex $v_1$ to $D$, $|N^{2}(D)|$ reduces by at least
$\lfloor\frac{\sigma_k}{2k}\rfloor+2$. When $G[N^{2}(D)]$ still has
an independent set $\{v_1,v_2,\cdots,v_k\}$ such that $ e(v_i,
N^{1}(D))\leq \lceil\frac{\sigma_k}{2k}\rceil-1$, we continue the
above manipulation, until $G[N^{2}(D)]$ has no independent set
$\{v_1,v_2,\cdots,v_k\}$ such that $ e(v_i, N^{1}(D))\leq
\lceil\frac{\sigma_k}{2k}\rceil-1$, where $i\in \{1,2,\cdots,k\}$.
Thus $D$ increases by at most
$2|N^{2}(D)|/(\lfloor\frac{\sigma_k}{2k}\rfloor+2)$. Hence $|D|\leq
3k\frac{n-|N^{2}(D)|-1}{\sigma_k+k}+3k-5+\frac{4k|N^{2}(D)|}{\sigma_k+4k}+1=
3k\frac{n-1}{\sigma_k+k}
-3k\frac{|N^{2}(D)|}{\sigma_k+k}+\frac{4k|N^{2}(D)|}{\sigma_k+4k}+3k-4<
3k\frac{n-1}{\sigma_k+k}+ \frac{k|N^{2}(D)|}{\sigma_k+4k}+3k-4
<4k\frac{n-1}{\sigma_k+k}+3k-4$. So $|D|
<4k\frac{n-1}{\sigma_k+k}+3k-4$.

Here, $ N^{2}(D)$ can be partitioned into two parts $N_1^{2}(D)$ and
$N_2^{2}(D)$, for $\forall v\in N_1^{2}(D), e(v,$\\$ N^{1}(D))\geq
\lceil\frac{\sigma_k}{2k}\rceil$, and $\forall v\in N_2^{2}(D), e(v,
N^{1}(D))\leq \lceil\frac{\sigma_k}{2k}\rceil-1$ and $\alpha
(G[N_2^{2}(D)])\leq k-1$, where $| N_1^{2}(D)|\geq 0,
|N_2^{2}(D)|\geq 0$. In the same way as before, we may arrive at
$|D|<4k\frac{n-1}{\sigma_k+k}+3k-4+2(k-1)=4k\frac{n-1}{\sigma_k+k}+5k-6$,
where $D$ is a connected two-way two-step dominating set, and for
each vertex $v\in N^{2}(D)$, $e(v, N^{1}(D))\geq
\lceil\frac{\sigma_k}{2k}\rceil$. So the claim  is true. \qed \\

\noindent{\bf Case 3.2.}  There does not exist any independent set
$\{v_1,v_2,\cdots,v_k\}$ in $G[N^{2}(D)]$ such that $e(v_i,
N^{1}(D))\leq \lceil\frac{\sigma_k}{2k}\rceil-1, i\in \{ 1,2,\cdots,
k\}$.

We partition $ N^{2}(D)$ into two parts $N_1^{2}(D)$ and
$N_2^{2}(D)$, for $\forall v\in N_1^{2}(D), e(v_i, N^{1}(D)) \geq
\lceil\frac{\sigma_k}{2k}\rceil$, and $\forall v\in N_2^{2}(D), e(v_i,
N^{1}(D))\leq \lceil\frac{\sigma_k}{2k}\rceil-1$ and $\alpha (G[N_2^{2}(D)])\leq
k-1$, where $| N_1^{2}(D)|\geq 0, | N_2^{2}(D)|\geq 0$. Similar to
the proof of the latter part of Case 3.1, we can get $|D|\leq
3k\frac{n-|N^{2}(D)|-1}{\sigma_k+k}+3k-5 +2(k-1)\leq
3k\frac{n-2}{\sigma_k+k}+5k-7$. The claim is also true. \qed \\

\noindent{\bf Claim 4.} $G$ has a connected
$\lceil\frac{\sigma_k}{1.9k}\rceil$-strong two-step dominating set
$D$ such that $|D|< \frac{38k(n-1)}{9(\sigma_k+k)}+5k-6$.\\

\noindent{\bf Proof.} We still look at the set $D$ in Claim 1 of the
proof for Theorem 6. If for each vertex $v\in N^{2}(D)$, $e(v,
N^{1}(D))\geq \lceil\frac{\sigma_k}{1.9k}\rceil$, then $D$ is
exactly the required dominating set. So we assume that there exists
a $v\in N^{2}(D)$ such that $e(v, N^{1}(D))\leq
\lceil\frac{\sigma_k}{1.9k}\rceil-1$.\\

\noindent{\bf Case 4.1.} There exists an independent set
$\{v_1,v_2,\cdots,v_k\}$ in $G[N^{2}(D)]$ such that $ e(v_i,
N^{1}(D)) \leq \lceil\frac{\sigma_k}{1.9k}\rceil-1, i\in \{
1,2,\cdots, k\}$.

Suppose that $d(v_1)\geq d(v_2)\geq\cdots\geq d(v_k)$. Let
$P=v_1v^{1}_1v^{0}_1$ be a shortest $v_1-D$ path, where $v_1$ has at
least
$\lceil\frac{\sigma_k}{k}\rceil-\lceil\frac{\sigma_k}{1.9k}\rceil+1\geq
\lfloor\frac{9\sigma_k}{19k}\rfloor+1$ neighbors in $N^{2}(D)$. When
we put vertex $v_1$ to $D$, $| N^{2}(D)|$ reduces by at least
$\lfloor\frac{9\sigma_k}{19k}\rfloor+2$. If $G[N^{2}(D)]$ still has
an independent set $\{v_1,v_2,\cdots,v_k\}$ such that $e(v_i,
N^{1}(D))\leq \lceil\frac{\sigma_k}{1.9k}\rceil-1$, we continue the
above manipulation, until $G[N^{2}(D)]$ has no independent set
$\{v_1,v_2,\cdots,v_k\}$ such that $e(v_i, N^{2}(D))\leq
\lceil\frac{\sigma_k}{1.9k}\rceil-1$, where $i\in\{1,2,\cdots,k\}$.
Thus $D$ increases by at most
$2|N^{2}(D)|/(\lfloor\frac{9\sigma_k}{19k}\rfloor+2)$. Hence
$|D|\leq
3k\frac{n-|N^{2}(D)|-1}{\sigma_k+k}+3k-5+\frac{38k|N^{2}(D)|}{9\sigma_k+38k}+1
<\frac{38k}{9}\frac{n-1}{\sigma_k+k}+3k-4$.\\
Here, $ N^{2}(D)$ can be partitioned into two parts $N_1^{2}(D)$ and
$N_2^{2}(D)$, for $\forall v\in N_1^{2}(D),$ $e(v, N^{1}(D))\geq
\lceil\frac{\sigma_k}{1.9k}\rceil$, and $\forall v\in N_2^{2}(D),
e(v, N^{1}(D))\leq \lceil\frac{\sigma_k}{1.9k}\rceil-1$ and $\alpha
(G[N_2^{2}(D)])\leq k-1$, where $| N_1^{2}(D)|\geq 0, |
N_2^{2}(D)|\geq 0$.  As
before, we can get $|D|<
\frac{38k}{9}\frac{n-1}{\sigma_k+k}+3k-4+2k-2$ where $D$ is a
$\lceil\frac{\sigma_k}{1.9k}\rceil$-strong two-step dominating set. So the claim
is true. \qed \\

\noindent{\bf Case 4.2.}  There does not exist any independent set
$\{v_1,v_2,\cdots,v_k\}$ in $G[N^{2}(D)]$ such that $ e(v_i,
N^{1}(D))\leq \lceil\frac{\sigma_k}{1.9k}\rceil-1, i\in \{
1,2,\cdots, k\}$.

It is also easy to check that $|D|\leq
3k\frac{n-|N^{2}(D)|-1}{\sigma_k+k}+3k-5 + 2(k-1)\leq
3k\frac{n-2}{\sigma_k+k}+5k-7$. The claim is true. \qed \\

By Lemma 2 and the definition of $rvc(G)$, we may assume that $G$
has less than $n(\sigma_k-k+1)+\frac{kn}{\sigma_k+k}$ edges. And by
Claim 3, we may first construct a $\lceil\frac{\sigma_k}{2k}\rceil$-strong
two-step dominating set $D$ with $|D|<
4k\frac{n-1}{\sigma_k+k}+5k-6$, and then we partition $N^{1}(D)$
into two parts $N_1^{1}(D)$ and $N_2^{1}(D)$, where $N_1^{1}(D)$ are
those vertices with at least $\frac{1}{2k^2}(\sigma_k+k)^{2}-1$
neighbors in $N^{2}(D)$. So we have $|N_1^{1}(D)|<
\frac{2k^{2}n}{\sigma_k+k}$. Let $N_2^{1}(D)$ are those vertices
which have at least one neighbor in $N_1^{1}(D)$, $N_2^{2}(D)=
N^{2}(D)\setminus N_1^{2}(D)$. Therefore, $e(N_2^{2}(D),
N_1^{1}(D))=\emptyset$.

Now we assign distinct colors to each vertex of $D\cup N_1^{1}(D)$,
then we color $N_1^{2}(D)$ only with 9 fresh colors so that each
vertex of $N_1^{2}(D)$ chooses its color randomly and independently
from all other vertices of $N_1^{2}(D)$. The vertices of $N^{2}(D)$
remain uncolored. We will show that the above coloring of $G$
results in a rainbow vertex-connection. It is obvious that any two
vertices of $D$ are connected by a rainbow path, as each vertex of
$D$ has distinct colors. Similarly, every pair of vertices of
$N_1^{1}(D)$ are connected by a rainbow path. For any two vertices
$u_1,u_2$ of $ N_1^{2}(D)$, there exist vertices $w_1,w_2\in D$ such
that $u_1w_1\in E(G), u_2w_2\in E(G)$, and $D$ has a $w_1-w_2$
rainbow path. Hence there exists a $u_1-u_2$ rainbow path in $G$. In
the same way, for any two vertices of $ N_2^{1}(D)$ we can find a
rainbow path connecting them. For $u\in D$ and $v\in N_1^{1}(D)$,
there is a vertex $w\in D$ such that $vw\in E(G)$, and $D$ has a
$w-u$ rainbow path. So $G$ has a $u-v$ rainbow path. In the same
way, for any two vertices coming from respectively any two sets of
$D, N_1^{1}(D),N_1^{2}(D),N_2^{1}(D)$ and $N_2^{2}(D)$, $G$ has a
rainbow path connecting them.

Now it suffices to prove that every two vertices of $N_2^{2}(D)$ are
connected by a rainbow path, that is, for $N_1^{2}(D)$, there exists
a coloring with 9 colors such that every vertex of $N_2^{2}(D)$ has
at least two neighbors in $N_1^{2}(D)$ colored differently. Let
$P_v$ be the event that all the neighbors of $v$ in $N_1^{2}(D)$ are
assigned at least two distinct colors. Now we will prove $Pr[P_v]>
0$ for each vertex $v \in N_2^{2}(D)$. As $D$ is a
$\lceil\frac{\sigma_k}{2k}\rceil$-strong two-step dominating set, we
can fix a set $X(v)\subset N_1^{2}(D)$ of neighbors of $v$ with
$|X(v)|=\lceil\frac{\sigma_k}{2k}\rceil$. Let $Q_v$ be the event
that all of the vertices in $X(v)$ receive the same color. Thus,
$Pr[Q_v]\leq 9^{-\lceil\frac{\sigma_k}{2k}\rceil+1}$. As each vertex
of $N_1^{2}(D)$ has less than $\frac{1}{2k^2}(\sigma_k+k)^{2}-1$
neighbors in $N_2^{2}(D)$, we have that the event $Q_v$ is
independent of all other events $Q_w$ for $v\neq w$ but at most
$(\frac{1}{2k^2}(\sigma_k+k)^{2}-2)\lceil\frac{\sigma_k}{2k}\rceil$
of them. Since $e\cdot
9^{-\lceil\frac{\sigma_k}{2k}\rceil+1}((\frac{1}{2k^2}
(\sigma_k+k)^{2}-2)\lceil\frac{\sigma_k}{2k}\rceil+1)< 1$ for all
$\sigma_k\geq 8k$, by the Lov\'{a}sz Local Lemma, we have $Pr[P_v]>
0$ for each $v \in N_2^{2}(D)$. Hence, we have proved that for
$N_1^{2}(D)$, there exists a coloring with 9 colors such that every
vertex of $N_2^{2}(D)$ has at least two neighbors in $N_1^{2}(D)$
colored differently.

We know that the total number of colors we used is at most
$|D|+|D_1|+9\leq
4k\frac{n-1}{\sigma_k+k}+5k-7+\frac{2k^{2}n}{\sigma_k+k}-1+9 <
\frac{(2k^{2}+4k)n}{\sigma_k+k}+5k+1$. Therefore, we arrive at
$rvc(G)\leq\frac{(2k^{2}+4k)n}{\sigma_k+k}+5k$ for all $\sigma_k\geq
8k$.

In the following we still make use of the above $G$, but we will use
Claim 4 to construct a $\lceil\frac{\sigma_k}{1.9k}\rceil$-strong
two-step dominating set $D$ with $|D|<
\frac{38k(n-1)}{9(\sigma_k+4k)}+5k-6$, and then we partition
$N^{1}(D)$ into two parts $N_1^{1}(D)$ and $N_1^{2}(D)$, where
$N_1^{1}(D)$ are those vertices with at least
$\frac{1}{2k^2}(\sigma_k+k)^{2}-1$ neighbors in $N^{2}(D)$. So we
have $|N_1^{1}(D)|< \frac{2k^{2}n}{\sigma_k+k}$. Let $N_2^{1}(D)$
consist of those vertices which have at least one neighbor in
$N_1^{1}(D)$, $N_2^{2}(D)= N^{2}(D)\setminus N_1^{2}(D)$.

Similar to the above coloring, we assign distinct colors to each
vertex of $D\cup  N_1^{1}(D)$, then we color $N_1^{2}(D)$ only with
9 fresh colors so that each vertex of $N_1^{2}(D)$ chooses its color
randomly and independently from all other vertices of $N_1^{2}(D)$.
The vertices of $N^{2}(D)$ remain uncolored. We will show that the
above coloring of $G$ results in a rainbow vertex-connection. We
only need to prove that every two vertices of $N_2^{2}(D)$ are
connected by a rainbow path. Let $P_v$ be the event that all the
neighbors of $v$ in $N_1^{2}(D)$ are assigned at least two distinct
colors. We will prove $Pr[P_v]> 0$ for each $v \in N_2^{2}(D)$. As
$D$ is a $\lceil\frac{\sigma_k}{1.9k}\rceil$-strong two-step
dominating set, we can fix a set $X(v)\subset N_1^{2}(D)$ of
neighbors of $v$ with $|X(v)|=\lceil\frac{\sigma_k}{1.9k}\rceil$.
Let $Q_v$ be the event that all of the vertices in $X(v)$ receive
the same color. Thus, $Pr[Q_v]\leq
9^{-\lceil\frac{\sigma_k}{1.9k}\rceil+1}$. As each vertex of
$N_1^{2}(D)$ has less than $\frac{1}{2k^2}(\sigma_k+k)^{2}-1$
neighbors in $N_2^{2}(D)$, we have that the event $Q_v$ is
independent of all other events $Q_w$ for $v\neq w$ but at most
$(\frac{1}{2k^2}(\sigma_k+k)^{2}-2)\lceil\frac{\sigma_k}{1.9k}\rceil$
of them. Since $e\cdot
9^{-\lceil\frac{\sigma_k}{1.9k}\rceil+1}((\frac{1}{2k^2}
(\sigma_k+k)^{2}-2)\lceil\frac{\sigma_k}{1.9k}\rceil+1)< 1$ for all
$\sigma_k\geq 7k+1$, by the Lov\'{a}sz Local Lemma, we have
$Pr[P_v]> 0$ for each $v \in N_2^{2}(D)$. Hence, we have proved that
for $N_1^{2}(D)$, there exists a coloring with 9 colors such that
every vertex of $N_2^{2}(D)$ has at least two neighbors in
$N_1^{2}(D)$ colored differently. And the total number of colors we
used is at most $|D|+|D_1|+9\leq
\frac{38k(n-1)}{9(\sigma_k+k)}+5k-7+\frac{2k^{2}n}{\sigma_k+k}-1+9 <
\frac{(2k^{2}+\frac{38}{9}k)n}{\sigma_k+k}+5k+1$. Therefore, we
arrive at $rvc(G)\leq\frac{(2k^{2}+\frac{38}{9}k)n}{\sigma_k+k}+5k$
for all $\sigma_k\geq 7k+1$.

Therefore, the proof of Theorem 9 is now complete. \qed

\end{document}